\documentclass[a4paper,11pt]{article}

\usepackage{latexsym}
\usepackage{amsmath,amsfonts,amssymb,amsthm}
\newtheorem{thm}{Theorem}[section]
\newtheorem{lem}[thm]{Lemma}

\newtheorem{cor}[thm]{Corollary}
\newtheorem{prop}[thm]{Proposition}

\newtheorem{rem}[thm]{Remark}

\input epsf

\title{On polytopes associated to factorisations of prime-powers}
\author{Roland Bacher}
\begin{document}
\maketitle

{\sl Abstract\footnote{ 
Math. class: 52B05, 52B11, 52B20}: 
We study polytopes associated to factorisations of prime powers. 
These polytopes have explicit
descriptions either in terms of their vertices or as
intersections of closed halfspaces associated to their facets.
We give formulae for their $f-$vectors.}


\section{Main results}\label{sectmainresults}
Polytopes have two dual descriptions: They can be given either as 
convex hulls of finite sets or as compact sets of the form
$\cap_{f \in\mathcal F} f^{-1}(\mathbb R_+)$ where 
$\mathcal F$ is a finite set of affine functions and where 
$f^{-1}(\mathbb R_+)$ denotes the closed half-space on which the
affine function $f$ is non-negative.

It is difficult to construct families of polytopes where both
descriptions are explicit. The aim of this paper is to study a
new family of such examples. These polytopes are associated to 
vector-factorisations 
of prime-powers where a {\it $d-$dimensional  vector-factorisation}
of a prime-power 
$p^e$ is an integral vector $(v_1,v_2,\dots,v_d)\in\mathbb N^d$ 
such that $p^e=v_1\cdot v_2\cdots v_d$.
Given a prime power $p^e\in\mathbb N$ and a 
natural integer $d\geq 1$, we denote by 
$\mathcal P(p^e,d)$ the convex hull of all 
$d-$dimensional vector-factorisations
$(v_1,v_2,\dots,v_d)\in\mathbb N^d$ of $p^e$.
The case $e=0$ yields the unique vector-factorisation $(1,1,\dots,1)$ 
and is without interest. 
For $d=2$ and $e\geq 2$ the polytope 
$\mathcal P(p^e,2)$ is a $2-$dimensional polygon with
vertices $(1,p^e),(p,p^{e-1}),\dots,(p^{e-1},p),(1,p^e)$.
For $e=1$, the polytope 
$\mathcal P(p,d)$ is a $(d-1)-$dimensional simplex with vertices
$(p,1,\dots,1),(1,p,1,\dots,1),\dots,(1,\dots,1,p)$.

The observation that the combinatorial properties of 
$\mathcal P(p^e,d)$ are independent of 
the prime $p$  in these examples is a general fact: 
The combinatorial properties of the polytope $\mathcal P(p^e,d)$
are independent always independent of the prime number $p$.
It is in fact possible to replace every occurence of $p$ by an arbitrary 
real constant which is strictly greater than $1$.
(The choice of a strictly positive real number which is strictly smaller than
$1$ leads to a combinatorially
equivalent polytope with the opposite orientation.) 
Let us also mention that the polytopes $\mathcal P(p^e,d)$
are invariant under permutations of coordinates.

In the sequel we suppose always $e\geq 2$ and $d\geq 2$. This ensures that
$\mathcal P(p^e,d)$ is $d-$dimensional. 

In order to state our first main result, the description of 
$\mathcal P(p^e,d)$ in terms of inequalities, we consider
for $\lambda\in\{1,\dots,\min(e,d-1)\}$ the set
$\mathcal R_\lambda(d,e)$ consisting of all integral vectors 
$\alpha=(\alpha_1,\dots,\alpha_d)\in\mathbb N^d$ such that 
$\min(\alpha_1,\dots,\alpha_d)=0$,
$$\max(\alpha_1,\dots,\alpha_d)\ d<e+\sum_{i=1}^d\alpha_i$$
and
$$e+\sum_{i=1}^d\alpha_i\equiv \lambda\pmod d\ .$$
We call $\mathcal R_\lambda(d,e)$ the set of {\it regular vectors
of type $\lambda$}. The union
$$\mathcal R(d,e)=\cup_{\lambda=1}^{\min(e,d-1)}\mathcal R_\lambda(d,e)$$
is the set of {\it regular vectors}. 

We denote by $\mu$ the
function on $\mathcal R(d,e)$ defined by the equality
$$\mu(\alpha)\ d+\lambda=e+\sum_{i=1}^d\alpha_i$$
where $\alpha=(\alpha_1,\dots,\alpha_d)$ is an element
of $\mathcal R_\lambda(d,e)$. The formula
$$\mu(\alpha)=\left\lfloor\left(e+\sum_{i=1}^d\alpha_i\right)/d\right\rfloor=
\left( e-\lambda+\sum_{i=1}^d \alpha_i\right)/d$$ 
shows that $\mu(\alpha)$ is a natural integer.

\begin{thm} \label{thminequalities}
Let $e\geq 2$ and $d\geq 2$ be two integers.

The polytope $\mathcal P(p^e,d)$ is defined by the inequalities
$$\sum_{i=1}^dx_i\leq (d-1)+p^e\ ,$$
$$x_i\geq 1,\ i=1,\dots,d\ ,$$
$$\sum_{i=1}^d p^{\alpha_i}x_i\geq\lambda p^{\mu(\alpha)+1}+
(d-\lambda)p^{\mu(\alpha)},\ \alpha=(\alpha_1,\dots,\alpha_d)\in
\mathcal R_\lambda(d,e)$$
for $\lambda$ in $\{1,\dots,\min(e,d-1)\}$.

This list of inequalities is minimal if $d\geq 3$. For $d=2$,
the minimal list is obtained by removing
the two inequalities $x_1\geq 1$ and $x_2\geq 1$.
\end{thm}

For $\lambda\in\{1,\dots,d-1\}$, we denote by 
$$\Delta(\lambda)=
\mathrm{Conv}((\epsilon_1,\dots,\epsilon_d)\in\{0,1\}^d,\ 
\sum_{i=1}^d \epsilon_i=\lambda)$$
the {\it $(d-1)-$dimensional hypersimplex of parameter 
$\lambda$} (see for example page 19 of \cite{Z}).

Faces of codimension $1$ (or $(d-1)-$dimensional faces) of a
$d-$dimensional polytope are called {\it facets}.
The following result describes all facets of $\mathcal P(p^e,d)$ 
in terms of their vertices.

\begin{thm} \label{thmfacets}
The set of all facets of $\mathcal P(p^e,d)$ (for $e\geq 2$
and $d\geq 3$) is given by the set of convex hulls of the following sets:

$$\{(p^e,1,\dots,1),(1,p^e,1,\dots,),\dots,(1,\dots,1,p^e)\},\ $$
$$\{(p^{\beta_1},\dots,p^{\beta_d})\in\mathbb N^d\ \vert \beta_i=0,\sum_{j=1}^d
\beta_j=e\},\ i=1,\dots,d$$
$$\{(p^{\beta_1+\epsilon_1},\dots,p^{\beta_d+\epsilon_d})\in\mathbb N^d\ \vert
\  (\epsilon_1,\dots,\epsilon_d)\in\Delta(\lambda)\}$$
with $\lambda=1,\dots,\min(e,d-1)$ and $(p^{\beta_1},\dots, p^{\beta_d})$
going through all $d-$dimensional vector-factorisations of $p^{e-\lambda}$.
\end{thm}

Recall that the $f-$vector of a $d-$dimensional polytope $\mathcal P$
counts the number $f_k$ of $k-$dimensional faces contained in $\mathcal P$.

The following result describes the coefficients of the $f-$vector of 
$\mathcal P(p^e,d)$.

\begin{thm} \label{thmfvector}
Let $e\geq 2$ and $d\geq 2$ be two integers.

The numbers $f_0,\dots,f_k$ with $f_k$ counting the number 
of $k-$dimensional faces of the polytope 
$\mathcal P(p^e,d)$ are given by the formulae
$$\begin{array}{l}
\displaystyle f_0={e+d-1\choose d-1}\ ,\\
\displaystyle f_1={d\choose 2}+{d\choose 2}{e+d-2\choose d-1}\ ,\\
\displaystyle f_k
={d+1\choose k+1}+{d\choose k+1}{e+d-1\choose d}-{d\choose k+1}
{e-k+d-1\choose d},\ 2\leq k<d\\
\displaystyle f_d=1\ .\end{array}$$
\end{thm}

The formula for the number $f_0$ of vertices is easy:
Identification of a vector-factorisation $(p^{\beta_1},p^{\beta_2},\dots,
p^{\beta_d})$ with the monomial $x_1^{\beta_1}x_2^{\beta_2}\cdots
x_d^{\beta_d}$ of $\mathbb Q[x_1,\dots,x_d]$ shows that 
$f_0$ is the dimension of the vector space spanned by homogeneous 
polynomials of degree $e$ in $d$ variables.

The plan of the paper is as follows:

Section \ref{sectgen} describes a few generalisations
of the polytopes $\mathcal P(p^e,d)$.

The next three sections are devoted to the proof of Theorem 
\ref{thminequalities} and Theorem \ref{thmfacets}.
The idea for proving Theorem \ref{thminequalities} is as follows: 
We show first that all inequalities of  
Theorem \ref{thminequalities} hold and that they correspond to facets 
of $\mathcal P(p^e,d)$.
Section \ref{sectregfacets} contains the details
for facets associated to regular vectors, called {\it 
regular facets}. Section \ref{sectexcfacets} is devoted 
to $d+1$ other obvious facets, called {\it exceptional facets}. 
The proofs contain explicit descriptions of all facets
and imply easily Theorem \ref{thmfacets} from  
Theorem \ref{thminequalities}. 
It remains to show that $\mathcal P(p^e,d)$ has no ``exotic'' facets
(neither regular nor exceptional).
This is achieved in Section \ref{sectallfacets}
by showing that any facet $f$ of the $(d-1)-$dimensional polytope
defined  
by a (regular or exceptional) facet $F$ of $\mathcal P$ is also contained 
in a second (regular or exceptional) facet $F'$ of $\mathcal P$.
This implies the completeness of the list of facets.

Finally, Section \ref{sectfvector} contains a proof of Theorem 
\ref{thmfvector}.

\section{Generalisations}\label{sectgen}

A straightforward generalisation of the polytope $\mathcal P(p^e,d)$
is obtained by considering the polytopes with vertices
given by vector-factorisations of an 
arbitrary natural integer $N$. The vertices of 
such a polytope $\mathcal P(N,d)$ are the 
$$\prod_{j=1}^h {e_j+d-1\choose d-1}$$
different $d-$dimensional vector-factorisations of 
$N=p_1^{e_1}\dots p_h^{e_h}$ where $p_1<p_2<\dots<p_h$ are all 
prime-divisors of $N$. A further generalisation is given by replacing 
$N$ with a monomial $\mathbf T=T_1^{e_1}\cdots T_h^{e_h}$ and by 
considering real evaluations $T_i=t_i>0$ of 
all $d-$dimensional vector-factorisations of $\mathbf T$
(defined in the obvious way). The combinatorial
type of these polytopes depends on these evaluations (or on the
primes $p_1,\dots,p_h$ involved in $N=p_1^{e_1}\cdots p_h^{e_h}$).
There should however exist a ``limit-type'' if the increasing sequence 
$p_1<p_2<\dots<p_h$ formed by all prime-divisors of $N$ grows extremely fast.

\begin{rem} One can also consider polytopes defined as the convex hull of 
vector-factorisations (of a given integer) subject to various 
restrictions. A perhaps interesting case is given by considering
only factorisations with decreasing coordinates. The number of
such decreasing $d-$dimensional vector-factorisations of $p^e$
equals the number of partitions of $e$ having at most $d$ 
parts if $N=p^e$ is the $e-$th power of a prime $p$.
\end{rem}

The polytopes $\mathcal P(N,d)$ have the following natural 
generalisation: Consider an integral symmetric matrix $A$ of size $d\times d$.
For a given natural number $N$, consider the set $\mathcal D_A(N)$ of
all integral diagonal
matrices $D$ of size $d\times d$ such that $A+D$ is positive definite 
and has determinant $N$. One can show that
$\mathcal D_A(N)$ is always finite.
Brunn-Minkowski's inequality for mixed volumes, see eg.
Theorem 6.2 in \cite{SY}, states that
$$\det\left(A+\sum_{D\in\mathcal D_A(N)}\lambda_D D\right)^{1/d}\geq
\sum_{D\in\mathcal D_A(N)}\lambda_D \det(A+D)^{1/d}=N^{1/d}$$
if $\sum_{D\in\mathcal D_A(N)}\lambda_D=1$ with $\lambda_D\geq 0$.
This inequality is strict except in
the obvious case where $\lambda_D$ is equal to $1$ for a unique 
matrix $D\in\mathcal D_A(N)$.
This implies that $\mathcal D_A(N)$ is the set of vertices of
the polytope $\mathcal P_A(N)$ defined as the convex hull of
$\mathcal D_A(N)$.

The polytope
$\mathcal P(N,d)$ discussed above
correspond to the case where $A$ is the zero matrix of size $d\times d$.

The choice of a Dynkin matrix of size $d\times d$ for a root system of
type $A$ and of $N=1$ leads to polytopes having ${2d\choose d}/(d+1)$
vertices, see for example the solution of problem (18) in 
\cite{CC} or \cite{LN}. 
These polytopes are different from the Stasheff polytopes
(or associahedra) since they are of dimension $d$ for $d\geq 3$.

Another natural and perhaps interesting choice is given by considering 
for $A$ (an integral multiple of) the all one matrix.

The determination of the number of vertices 
of $\mathcal P_A(N)$ (or even of $\mathcal P_A(1)$) is perhaps 
a non-trivial problem, say, for $A$ the adjacency matrix of a 
connected finite simple graph.

\section{Regular facets}\label{sectregfacets}

We show first that all inequalities of Theorem \ref{thminequalities}
associated to regular vectors in $\mathcal R_\lambda(d,e)$ 
are satisfied on $\mathcal P(p^e,d)$.

We show then that every such inequality is 
sharp on a subset of vertices in $\mathcal P(p^e,d)$
spanning a polytope affinely equivalent to a $(d-1)-$dimensional
hypersimplex. All these inequalities define thus facets and are
necessary. We call a facet $F_\alpha$ associated to such a regular
vector $\alpha\in\mathcal R_\lambda(p^e,d)$ a {\it regular facet}.

\begin{prop} \label{propineqregular}
Let $\alpha\in \mathcal R_\lambda(d,e)$ 
be a regular vector. Setting
$\mu=\mu(\alpha)$, we have 
$$\sum_{i=1}^d p^{\alpha_i}x_i\geq \lambda p^{\mu+1}+(d-\lambda)p^\mu$$
for every element $(x_1,\dots,x_d)$ of $\mathcal P(p^e,d)$.
\end{prop}

{\bf Proof} Given $\alpha=(\alpha_1,\dots,\alpha_d)$ in $\mathcal 
A_\lambda=\mathcal R_\lambda(d,e)$, 
we denote by $l_\alpha$ the linear form defined by 
$$l_\alpha(x_1,\dots,x_d)=\sum_{i=1}^dp^{\alpha_i}x_i\ .$$
We have to show that $l_\alpha(x)\geq \lambda p^{\mu+1}+(d-\lambda)p^\mu$
for $x$ in $\mathcal P=\mathcal P(p^e,d)$ and $\mu=\mu(\alpha)$.
Since $l_\alpha$ is linear, it is enough to establish the inequality 
for all vertices of $\mathcal P$.

Let $v=(p^{\beta_1},\dots,p^{\beta_d})$ be a vertex of $\mathcal P$
realising the minimum $l_\alpha(v)=\min l_\alpha(\mathcal P)$. Set 
$$\begin{array}{l}
\displaystyle a=\min(\alpha_1+\beta_1,\dots,\alpha_d+\beta_d)\ ,\\
\displaystyle A=\max(\alpha_1+\beta_1,\dots,\alpha_d+\beta_d)\ .\end{array}$$
Choose indices $i,j$ in $\{1,\dots,d\}$ such that
$a=\alpha_i+\beta_i$ and $A=\alpha_j+\beta_j$.

If $a=A$, the equality 
$$dA=\sum_{i=1}^d\alpha_i+\beta_i=e+\sum_{i=1}^d\alpha_i\equiv \lambda
\pmod d$$
shows $\lambda=0$ in contradiction with $\lambda\in\{1,\dots,d-1\}$.

We claim next that $\beta_j\geq 1$. Indeed, we have otherwise
$A=\alpha_j$ and
$$e+\sum_{k=1}^d\alpha_k=\sum_{k=1}^d(\alpha_k+\beta_k)\leq A d= 
\max(\alpha_1,\dots,\alpha_d)\ d$$
in contradiction with the inequality 
$\max(\alpha_1,\dots,\alpha_d)\ d<e+\sum_{k=1}^d \alpha_k$ satisfied by
$\alpha\in\mathcal R_\lambda$.

We consider now the vertex 
$$\tilde v=(p^{\tilde \beta_1},\dots,p^{\tilde \beta_d})$$
where $\tilde\beta_k=\beta_k$ if $k\not\in \{i,j\},\ \tilde 
\beta_i=\beta_i+1=a+1$ and $\tilde \beta_j=\beta_j-1=A-1$.
We have
$$l_\alpha(v)-l_\alpha(\tilde v)=\sum_{k=1}^d p^{\alpha_k+\beta_k}-
\sum_{k=1}^d p^{\alpha_k+\tilde\beta_k}=p^A+p^a-(p^{A-1}+p^{a+1})\ .$$
Since $p>1$ and $A>a$ we have
$$p^A+p^a-(p^{A-1}+p^{a+1})=(p^{A-1}-p^a)(p-1)\geq 0\ .$$
This shows that $A-1=a$ by minimality of $l_\alpha(v)$.

In order to compute the value of $l_\alpha(v)$, we use 
the equalities
$$\sum_{i=1}^d\alpha_i+\beta_i=e+\sum_{i=1}^d\alpha_i=\lambda+\mu d\ .$$
Since the vector 
$w=(\alpha_1+\beta_1,\dots,\alpha_d+\beta_d)$ has all its coefficients 
in $\{a,a+1=A\}$, we get $a=\mu$ and $w$ takes the value $\mu$ with 
multiplicity $d-\lambda$ and $\mu+1$ with multiplicity $\lambda$.
This shows $l_\alpha(v)=\lambda p^{\mu+1}+(d-\lambda)p^\mu$.
\hfill$\Box$

\begin{prop} \label{propbij}
The map $(\alpha_1,\dots,\alpha_d)\longmapsto
(p^{\mu-\alpha_1},\dots,p^{\mu-\alpha_d})$
(where $\mu=\lfloor (e+\sum_{i=1}^d \alpha_i)/d\rfloor=
(e-\lambda+\sum_{i=1}^d\alpha_i)/d$) is a one-to-one map from
the set $\mathcal R_\lambda(d,e)$ of regular vectors of type $\lambda$
onto the set of $d-$dimensional vector-factorisations of $p^{e-\lambda}$.
The inverse map is given by $(p^{\beta_1},\dots,p^{\beta_d})
\longmapsto (B-\beta_1,\dots,B-\beta_d)$
where $B=\max(\beta_1,\dots,\beta_d)$.
\end{prop}

{\bf Proof} 
We define $\mathcal F_\lambda=\mathcal F_\lambda(d,e)$ as the 
finite set of all integral vectors 
$(\beta_1,\dots,\beta_d)\in\mathbb N^d$ such that 
$\sum_{i=1}^d\beta_i=e-\lambda$. The map
$$\mathcal F_\lambda\ni(\beta_1,\dots,\beta_d)\longmapsto
(p^{\beta_1},\dots,p^{\beta_d})$$
yields a bijection between $\mathcal F_\lambda$ and
the set of $d-$dimensional vector-factorisations of $p^{e-\lambda}$.

We show first the inclusions $\varphi(\mathcal R_\lambda)\subset
\mathcal F_\lambda$ and $\psi(\mathcal F_\lambda)\subset
\mathcal R_\lambda$ where
$$\varphi(\alpha_1,\dots,\alpha_d)=(\mu-\alpha_1,\dots,\mu-\alpha_d)$$
with $\mu=\left(e-\lambda+\sum_{i=1}^d\alpha_i\right)/d$ and 
$$\psi(\beta_1,\dots,\beta_d)=(B-\beta_1,\dots,B-\beta_d)$$
with $(\beta_1,\dots,\beta_d)\in\mathbb N^d$ such that 
$\sum_{i=1}^d \beta_i=e-\lambda$ and $B=\max(\beta_1,\dots,\beta_d)$.

The inclusion $\varphi(\mathcal R_\lambda)\subset
\mathcal F_\lambda$ follows from 
$$\max(\alpha_1,\dots,\alpha_d)\leq \lfloor(e+\sum_{i=1}^d\alpha_i)/d\rfloor
=\mu=(e-\lambda+\sum_{i=1}^d\alpha_i)/d$$
showing $\varphi(\alpha_1,\dots,\alpha_d)\in\mathbb N^d$ and from 
$$\sum_{i=1}^d(\mu-\alpha_i)=\mu d-\sum_{i=1}^d\alpha_i=e-\lambda\ .$$

Consider now $(\beta_1,\dots,\beta_d)\in \mathcal F_\lambda$.
We have $(B-\beta_1,\dots,B-\beta_d)\in \mathbb N^d$ and 
$\min(B-\beta_1,\dots,B-\beta_d)=0$ where $B=\max(\beta_1,\dots,\beta_d)$.
We have moreover the inequalities
$$e+\sum_{i=1}^d(B-\beta_i)=\lambda+Bd>Bd\geq 
\max(B-\beta_1,\dots,B-\beta_d)\ d$$
since $\lambda>0$. Finally, the computation
$$e+\sum_{i=1}^d(B-\beta_i)=\lambda+Bd\equiv \lambda\pmod d$$
proves the inclusion of $\psi(\beta_1,\dots,\beta_d)=(B-\beta_1,\dots,
B-\beta_d)$ in $\mathcal R_\lambda$.

The computation of 
$$\mu=\left(e-\lambda+\sum_{i=1}(B-\beta_i)\right)/d=B$$
shows that $\varphi\circ \psi$ is the identity map of $\mathcal F_\lambda$.

Consider $(\alpha_1,\dots,\alpha_d)\in \mathcal R_\lambda$.
Since $\min(\alpha_1,\dots,\alpha_d)=0$, we have 
$\max(\mu-\alpha_1,\dots,\mu-\alpha_d)=\mu$. This implies that 
$$\psi\circ\varphi(\alpha_1,\dots,\alpha_d)=(\mu-(\mu-\alpha_1),\dots,
\mu-(\mu-\alpha_d))$$ 
is the identity map of the set $\mathcal R_\lambda$.
\hfill$\Box$

\begin{prop} \label{propaffhypersimplex}
(i) Consider an integral vector $(\beta_1,\dots,\beta_d)\in\mathbb N^d$,
an integer $\lambda\in\{1,\dots,d-1\}$ and a prime $p$.
The affine isomorphism
$$(x_1,\dots,x_d)\longmapsto \left(\frac{
x_1-p^{\beta_1}}{p^{\beta_1+1}-p^{\beta_1}},\dots,
\frac{x_d-p^{\beta_d}}{p^{\beta_d+1}-p^{\beta_d}}\right)$$
of $\mathbb R^d$ induces a one-to-one map between the set 
$$S=\{(p^{\beta_1+\epsilon_1},\dots,p^{\beta_d+\epsilon_d})\in
\mathbb N^d\ \vert\
(\epsilon_1,\dots,\epsilon_d)
\in\Delta(\lambda)\}$$
and the set 
$$\left\{(\epsilon_1,\dots,\epsilon_d)
\in\{0,1\}^d,\ \sum_{i=1}^d \epsilon_i=\lambda\right\}$$
of vertices of the $(d-1)-$dimensional hypersimplex 
$\Delta(\lambda)$.

\ \ (ii) The facets of the convex hull of $S$
are of the form $x_i=p^{\beta_i+\epsilon}$ for $i=1,\dots,d$ and $\epsilon\in
\{0,1\}$ except if $\lambda=1$ or $\lambda=d-1$ where all facets are of the
form $x_i=p^{\beta_i}$ respectively $x_i=p^{\beta_i+1}$.
\end{prop}

{\bf Proof} We leave the easy proof of assertion (i) to the reader.

Assertion (ii) follows from the peculiar form of the affine isomorphism
introduced in assertion (i) and from the observation 
that facets of $\Delta(\lambda)$
are given as intersections of the hyperplane defined by the equation
$\sum_{i=1}^d x_i=\lambda$ with one of the $2d$ facets of the $d-$dimensional
cube $[0,1]^d$.\hfill $\Box$

\begin{cor} \label{corfacet} The inequality 
$$\sum_{i=1}^dp^{\alpha_i} x_i\geq\lambda p^{\mu+1}+(d-\lambda)p\mu$$
associated to a regular vector $(\alpha_1,\dots,\alpha_d)$ in 
$\mathcal R_\lambda(d,e)$
is sharp on the set
$$S_\alpha=\{(p^{\mu-\alpha_1+\epsilon_1},\dots,p^{\mu-\alpha_d
+\epsilon_d})\ 
\vert\ (\epsilon_1,\dots,\epsilon_d)\in \Delta(\lambda)\}$$
of vertices of $\mathcal P(p^e,d)$.
The convex hull of $S_\alpha$ is a facet of $\mathcal P(p^e,d)$
which is affinely equivalent to the
$(d-1)-$dimensional hypersimplex $\Delta(\lambda)$.
\end{cor}

{\bf Proof} The obvious equalities
$$\begin{array}{l}
\displaystyle
\mu=\min(\mu-\alpha_1+\epsilon_1+\alpha_1,\dots,\mu-\alpha_d+\epsilon_d
+\alpha_d)\\
\displaystyle  
\mu+1=\max(\mu-\alpha_1+\epsilon_1+\alpha_1,\dots,\mu-\alpha_d+\epsilon_d
+\alpha_d)\end{array}$$
and the definition of $\Delta(\lambda)$ show that $S_\alpha$ 
consists exactly of all vertices of $\mathcal P(p^e,d)$ such that $A=a+1$ 
with $a,A$ as in the proof of Proposition \ref{propineqregular}.
The arguments of the proof of Proposition \ref{propineqregular}
imply thus that $S_\alpha$ is the subset of vertices of
$\mathcal P(p^e,d)$ on which the linear form
$$x=(x_1,\dots,x_d)\longmapsto l_\alpha(x)=\sum_{i=1}^d p^{\alpha_i}x_i$$
is minimal. This shows that the convex hull $F_\alpha$ of the set $S_\alpha$ 
defines a $k-$dimensional face of $\mathcal P(p^e,d)$ for some integer
$k\in\{0,\dots,d-1\}$. Assertion (i) of Proposition 
\ref{propaffhypersimplex} implies now that $F_\alpha$ is affinely
equivalent to a hypersimplex $\Delta(\lambda)$ of 
dimension $d-1$. In particular, the convex hull $F_\alpha$ of 
$S_\alpha$ is a facet of $\mathcal P(p^e,d)$.
\hfill$\Box$

\section{Exceptional facets}\label{sectexcfacets}

We leave it the reader to check that we have 
$$\sum_{i=1}^dx_i\leq d-1+p^e$$
for $(x_1,\dots,x_d)\in\mathcal P(p^e,d)$. The details are straightforward
and involve computations similar to those used for 
proving  Proposition \ref{propineqregular}.
Equality holds for the elements of the exceptional facet $F_\infty$
given by the $(d-1)-$dimensional simplex with vertices 
$(p^e,1,\dots,1),\dots,(1,\dots,1,p^e)$.

The inequalities $x_i\geq 1,\ i=1,\dots,d$ hold obviously for 
$(x_1,\dots,x_d)\in
\mathcal P(p^e,d)$. For $d\geq 3$,
these inequalities define $d$ exceptional facets $F_1,\dots,F_d$
which are all affinely equivalent to the $(d-1)-$dimensional
polytope $\mathcal P(p^e,d-1)$.

\begin{rem}
The $d+1$ inequalities associated to exceptional facets define
a $d-$dimensional simplex with vertices 
$(1,1,\dots,1,1),(p^e,1,\dots,1),\dots,(1,\dots,1,p^e)$.
This simplex contains $\mathcal P(p^e,d)$.
\end{rem}

\section{Proof of Theorem \ref{thminequalities} and \ref{thmfacets}}
\label{sectallfacets}

The main tool for proving Theorem \ref{thminequalities} 
is the following obvious and well-known result.

\begin{prop} \label{propcriterefaces} Let $\mathcal F$ be a non-empty 
set of facets of a polytope $\mathcal P$.
The set $\mathcal F$ contains all facets of $\mathcal P$
if and only if for every element $F\in \mathcal F$ and for every facet
$f$ of $F$, there exists a distinct element $F'\not=F$ in $\mathcal F$
such that $f$ is also a facet of $F'$.
\end{prop}

{\bf Proof} Call two facets of a $d-$dimensional polytope $\mathcal P$
{\it adjacent} if they intersect in a common $(d-2)-$face of $\mathcal P$.
Consider the graph with vertices formed by all facets of $\mathcal P$
and edges given by adjacent pairs of facets. This graph is connected
and its edges are in bijection with $(d-2)-$faces of $\mathcal P$
since the intersection of three distinct facets is of dimension $\leq d-3$.
Proposition \ref{propcriterefaces}
boils now down to the trivial observation that a non-empty 
subset $\mathcal V'$ of vertices of a connected graph $\Gamma$ coincides with 
the set of vertices of $\Gamma$ if and only if for every vertex $v$ of
$\mathcal V'$, the set $\mathcal V'$ contains also all vertices of $\Gamma$ 
which are adjacent to $v$.
\hfill $\Box$

Given a subset $S$ of vertices of $\mathcal P(p^e,d)$,
we consider
$$m_i=\min_{(p^{\beta_1},\dots,p^{\beta_d})\in S}(\beta_i)$$
and 
$$M_i=\max_{(p^{\beta_1},\dots,p^{\beta_d})\in S}(\beta_i)\ .$$
We have thus 
$$m_i\leq \beta_i\leq M_i$$
for every element $(p^{\beta_1},\dots,p^{\beta_d})$
of $S$ and these inequalities are sharp.
We set $m(S)=(m_1,\dots,m_d)$ and $M(S)=(M_1,\dots,M_d)$.

An important ingredient of all proofs is the following result.

\begin{lem} \label{lemregsubset}
Let $S$ be a subset of vertices of $\mathcal P(p^e,d)$ 
such that $M-m\in\{0,1\}^d$ 
where $m=m(S)$ and $M=M(S)$ are as above.
Then $S$ is contained in the set of vertices of a 
regular facet of $\mathcal P(p^e,d)$.
\end{lem}

{\bf Proof} Set $\lambda=e-\sum_{i=1}^dm_i$.
Suppose first $\lambda=0$. This implies that $S$ is reduced to a 
unique element $v=(p^{m_1},\dots,p^{m_d})$. Choose two distinct 
indices $i,j$ in $\{1,\dots,d\}$  
such that $m_i>0$ in order to construct  
the element $\tilde v=(p^{\tilde m_1},\dots,p^{\tilde m_d})$ 
where $\tilde m_k=m_k$ if $k\not=i,j$, $\tilde m_i=m_i-1,\tilde m_j=m_j+1$.
The set $\tilde S=\{v,\tilde v\}$ contains $S$ and
satisfies the conditions of Lemma \ref{lemregsubset}
with $\lambda=1$. We may thus assume $\lambda\geq 1$.

The obvious identity $\prod_{i=1}^d p^{\beta_i}=p^e$
shows the equality $\sum_{i=1}^d \beta_i=\lambda+\sum_{i=1}^d
m_i$ for every element $(p^{\beta_1},\dots,p^{\beta_d})$ of $S$.
The inclusion $M-m\in\{0,1\}^d$ shows 
$m_i\leq \beta_i\leq M_i\leq m_i+1$ and implies $\lambda\leq d$.
Moreover, if $\lambda=d$ then $\beta_i=m_i+1$
for every element $(p^{\beta_1},\dots,p^{\beta_d})$ of $S$.
This shows that $S$ is reduced to the unique element
$(p^{m_1+1},\dots,p^{m_d+1})$ and contradicts the definition 
of $m=(m_1,\dots,m_d)$. We have thus $\lambda\in\{1,\dots,d-1\}$.

Up to enlarging the set $S$, we can assume
$$S=\{(p^{m_1+\epsilon_1},\dots,p^{m_d+\epsilon_d})\in\mathbb N^d\ \vert\  
(\epsilon_1,\dots,\epsilon_d)\in \Delta(\lambda)\}$$
for some integer $\lambda\in\{1,\dots,d-1\}$.

Consider the function $\psi$ defined by $\psi(m_1,\dots,m_d)=
(\mu-m_1,\dots,\mu-m_d)$ where $\mu=\max(m_1,\dots,m_d)$, 
see the proof of Proposition \ref{propbij}. 
Proposition \ref{propbij} shows that we have
$\psi(m_1,\dots,m_d)=(\mu-m_1,\dots
\mu-m_d)\in\mathcal R_\lambda(p^e,d)$. 
Corollary \ref{corfacet} shows now that $\mathcal S$ is 
the vertex set of the regular facet defined by the regular vector
$\psi(m_1,\dots,m_d)=(\mu-m_1,\dots,\mu-m_d)$ of 
$\mathcal R_\lambda(p^e,d)$.
\hfill $\Box$

{\bf Proof of Theorem \ref{thminequalities}} 
Corollary \ref{corfacet} and Section \ref{sectexcfacets}
show that all inequalities
of Theorem \ref{thminequalities} are satisfied and necessary if 
$d\geq 3$. We have thus to show that every facet of $\mathcal P(p^e,d)$ 
is either in $\{F_\infty,F_1,F_2,\dots,F_d\}$ or is among the set
$\{F_\alpha\}_{\alpha\in\mathcal R}$ of
regular facets indexed by the set $\mathcal R=\cup_{\lambda=1}
^{\min(e,d-1)}\mathcal R_\lambda(d,e)$ of regular vectors.
We show this using Proposition \ref{propcriterefaces}
with respect to the set of facets 
$$\mathcal F=
F_\infty\cup\bigcup_{i=1}^d F_i\cup\bigcup_{\alpha\in\mathcal R} F_\alpha
\ .$$

We consider first the exceptional facet $F_\infty$. 
A facet $f$ of $F_\infty$ is defined by 
an additional equality $x_i=1$ for some $i\in\{1,\dots,d\}$
and we have thus $f\in F_\infty\cap F_i$ where $F_i$ is the 
exceptional facet defined by $x_i=1$.

Consider next an exceptional facet $F_i$ defined by $x_i=1$ for some
$i\in\{1,\dots,d\}$. Such a facet $F_i$ coincides with the polytope
$\mathcal P'=\mathcal P(p^e,d-1)$ of all $(d-1)-$dimensional
vector-factorisations of $p^e$.  
Facets of $F_i$ are thus in bijection with facets of $\mathcal P'$.
Using induction on $d$ (the initial case $d=2$ is easy), we know
thus complete list of facets of $F_i$. 
Consider first a facet $f$ of $F_i$ corresponding to an ordinary facet
of $\mathcal P'$. Its vertices satisfy the conditions of Lemma 
\ref{lemregsubset} and are thus also contained in a regular facet of 
$\mathcal P$.
A facet $f$ of $F_i$ corresponding to the exceptional facet $F'_\infty$
of $\mathcal P'$ is also contained in the exceptional facet 
$F_\infty$ of $\mathcal P$.
All other exceptional facets of $\mathcal P'$ 
are given by $x_i=x_j=1$ for 
some $j\not=i$ and are thus contained in $F_i\cap F_j$.
(Remark that the last case does never arise for $d=3$.)
 
We consider now a regular facet $F$ of type $\lambda$
with vertices 
$$\{(p^{\beta_1+\epsilon_1},\dots,p^{\beta_d+\epsilon_d})\ \vert \ 
(\epsilon_1,\dots,\epsilon_d)\in \Delta(\lambda)\}\  .$$
Since $F$
is affinely equivalent after multiplication by a diagonal matrix and 
a translation
to the hypersimplex $\Delta_{d-1}(\lambda)$,
a facet $f$ of $F$ is given by an additional equality $x_i=c$
with $i$ in $\{1,\dots,d\}$ and $c$ in $\{p^{\beta_i},
p^{\beta_i+1}\}$.

If $c=1$ then $f$ is also contained in the exceptional facet $F_i$.

If $c=p^{\beta_i}>1$ and $\lambda<d-1$ then $f$ belongs also to the 
regular facet 
of type $\lambda+1$ with vertices
$$\{(t^{\beta_1+\epsilon_1},\dots,p^{\beta_{i-1}+\epsilon_{i-1}},p^{
\beta_i-1+\epsilon_i},p^{\beta_{i+1}+\epsilon_{i+1}},\dots,p^{\beta_d+
\epsilon_d})\ \vert\ (\epsilon_1,\dots,\epsilon_d)\in \Delta(\lambda+1)\}\ .$$
The case $\gamma_i=\beta_i$ and $\lambda=d-1$ implies that 
the set of vertices of $f$ is reduced to 
$(p^{\beta_1+1},\dots,p^{\beta_{i-1}+1},p^{\beta_i},p^{\beta_{i+1}+1},
\dots,p^{\beta_d+1})$. This is impossible for $d\geq 3$.

If $c=p^{\beta_i+1}$ and $\lambda>1$ then $f$ belongs also to the 
regular facet of type $\lambda-1$ with vertices
$$\{(t^{\beta_1+\epsilon_1},\dots,p^{\beta_{i-1}+\epsilon_{i-1}},p^{
\beta_i+1+\epsilon_i},p^{\beta_{i+1}+\epsilon_{i+1}},\dots,p^{\beta_d+
\epsilon_d})\ \vert\ (\epsilon_1,\dots,\epsilon_d)\in \Delta(\lambda-1)\}\ .$$
In the case $c=p^{\beta_i+1}$ and $\lambda=1$ the set of vertices 
of $f$ is reduced to 
$(p^{\beta_1},\dots,p^{\beta_{i-1}},p^{\beta_i+1},p^{\beta_{i+1}},
\dots,p^{\beta_d})$ and this is impossible for $d\geq 3$.

Proposition \ref{propcriterefaces}
shows that $\mathcal F$ is the complete list of facets for 
$\mathcal P$. This ends the proof of Theorem \ref{thminequalities}.
\hfill$\Box$

{\bf Proof of Theorem \ref{thmfacets}} Theorem \ref{thmfacets}
 follows easily from Theorem \ref{thminequalities}
and from the explicit descriptions of regular facets given by
Proposition \ref{propbij} and Corollary \ref{corfacet}.

\section{Proof of Theorem \ref{thmfvector}}\label{sectfvector}

Theorem \ref{thmfvector} is easily checked in the case $d=2$
where $\mathcal P(p^e,2)$ is the polygon defined by the $e+1$ 
vertices $(p^e,1),(p^{e-1},p),\dots,(p,p^{e-1}),(1,p^e)$.

We suppose henceforth $e\geq 2$ and $d\geq 3$ and we
consider  $\mathcal P=\mathcal P(p^e,d)$ (where $p$ is a prime). 

The formula for the number $f_0$ of vertices of $\mathcal P$ certainly 
holds. Indeed, $f_0$ is equal to the number of $d-$dimensional
vector-factorisations of $p^e$. Such vector-factorisations are in 
one-to-one correspondence with homogeneous monomials of degree $e$
in $d$ commuting variables. The vector space spanned by these monomials
is of dimension ${e+d-1\choose d-1}$.

We call a $k-$dimensional face of $\mathcal P$ 
{\it regular} if it is contained in a
regular facet of $\mathcal P$. 
A $k-$dimensional face of $\mathcal P$ is {\it exceptional} otherwise. 

A $1-$dimensional face is exceptional if and only if it 
is contained in the exceptional facet $F_\infty$. There are thus
${d\choose 2}$ exceptional $1-$dimensional faces.

For $k\geq 2$, a $k-$dimensional face $f$ which is exceptional is
either contained in the exceptional facet $F_\infty$ and there are 
${d\choose k+1}$ such faces, or it is the intersection of 
$d-k$ distinct exceptional facets in $\{F_1,\dots,F_d\}$. 
For $k$ in $\{2,\dots,d-1\}$, the polytope $\mathcal P$ contains thus 
$${d\choose k+1}+{d\choose d-k}={d+1\choose k+1}$$
exceptional $k-$dimensional faces.

A regular face $f$ of dimension $k\geq 1$ has
a type $\lambda=e-\sum_{i=1}^d m_i\in \{1,\dots,k\}$ 
where $m(S)=(m_1,\dots,m_d)$ is associated to the vertex set 
$S$ of $f$ as in Lemma \ref{lemregsubset}. The face $f$ is then defined  
by the support consisting of the $k+1$ non-zero coordinates of
$M(S)-m(s)\in\{0,1\}^d$ and by the regular facet with vertices
$$\{(p^{m_1+\epsilon_1},\dots,p^{m_d+\epsilon_d})\ \vert\ 
(\epsilon_1,\dots,\epsilon_d)\in\Delta(\lambda)\}\ .$$
The number of regular $k-$dimensional faces contained in $\mathcal P$
is thus given by
$${d\choose k+1}\sum_{\lambda=1}^{\min(k,e)}\sharp(A_\lambda)=
{d\choose k+1}\sum_{\lambda=1}^{\min(k,e)}{e-\lambda+d-1\choose d-1}$$
where the first factor corresponds to the choice of a support for 
$M(S)-m(S)$, the sum corresponds to all possibilities for $\lambda$
and the factor $\sharp(A_\lambda)={e-\lambda+d-1\choose d-1}$
corresponds to all possibilities for the \lq\lq minimal'' regular
facet with vertices 
$$\{(p^{m_1+\epsilon_1},\dots,p^{m_d+\epsilon_d})\ \vert\ 
(\epsilon_1,\dots,\epsilon_d)\in\Delta(\lambda)\}$$
which contains $f$.

Iterated application of the identity 
${a-1\choose b-1}+{a-1\choose b}={a\choose b}$ shows 
$$\sum_{\lambda=1}^{\min(k,e)}{e-\lambda+d-1\choose d-1}={e-
1+d\choose d}-{e-\min(e,k)+d-1\choose d}\ .$$ 
This yields the closed expression
$$f_k={d+1\choose k+1}+{d\choose k+1}{e+d-1\choose d}-{d\choose k+1}
{e-k+d-1\choose d}$$
for $f_2,\dots,f_{d-1}$ and ends the proof of Theorem 
\ref{thmfvector}.\hfill$\Box$

\begin{rem} The proof of Theorem \ref{thmfvector} contains the 
detailled description in terms of vertices of all faces of 
$\mathcal P$. It is thus easy to work out the
face-lattice of $\mathcal P(p^e,d)$.
\end{rem}


I thank F. Mouton for helpful comments.

\end{document}